\documentclass[12pt]{amsart}
\usepackage{a4wide}
\usepackage{amssymb}
\newtheorem*{theorem*}{Theorem}

\newtheorem*{proposition*}{Proposition}
\newtheorem*{corollary}{Corollary}
\newtheorem{lemma}{Lemma}

\newtheorem*{mlemma*}{Main Lemma}
\newtheorem*{lemma*}{Lemma}

\newtheorem*{claim*}{Claim}
\theoremstyle{remark}

\newcommand{\nc}{\newcommand}
\newcommand{\rnc}{\renewcommand}

\newcommand{\C}{{\mathbb C}}

\newcommand{\R}{{\mathbb R}}

\nc{\ol}{\overline}
\nc{\bea}{\begin{eqnarray}}
\nc{\eea}{\end{eqnarray}}
\nc{\beqa}{\begin{eqnarray*}} \nc{\bega}{\begin{gather*}}
\nc{\Sd}{\stackrel{\cdot}{+}} \nc{\eeqa}{\end{eqnarray*}}
\nc{\eega}{\end{gather*}} \nc{\nn}{\nonumber}
\nc{\lmt}{\longmapsto} \nc{\eps}{\varepsilon}
\nc{\Llra}{\Longleftrightarrow} \nc{\LRA}{\Leftrightarrow}
\nc{\lra}{\longrightarrow} \nc{\Lra}{\Longrightarrow}
\nc{\Lla}{\Longleftarrow} \nc{\lmto}{\longmapsto}
\nc{\vt}{\vartheta} \nc{\vs}{\psi} \nc{\sn}{\sigma_n}
\nc{\Hi}{H^{\infty}} \nc{\Li}{L^{\infty}} \nc{\vp}{\varphi}
\nc{\Tp}{T_{\vp}} \nc{\Tpb}{T_{\overline{\vp}}}
\nc{\Tpn}{T_{\vp_0}} \nc{\Tpnb}{T_{\overline{\vp}_0}}
\nc{\gp}{g^{2/p}} \nc{\Kpw}{K^p_I(|g|^2)} \nc{\Kqw}{K^q_I(|g|^2)}
\nc{\Kp}{K^p_I} \nc{\kep}{\operatorname{ker} T_{\vp}}
\nc{\Lpw}{L^p(|g|^2)} \nc{\Lppw}{L^{p'}(|g|^2)}
\nc{\Lqw}{L^q(|g|^2)} \nc{\Ktw}{K^2_I(|g|^2)} \nc{\Kt}{K^2_I}
\nc{\Ltw}{L^2(|g|^2)} \nc{\kl}{k_{\lambda}}
\nc{\klb}{k^b_{\lambda}} \nc{\klIb}{k^{Ib}_{\lambda}}
\nc{\klI}{k_{\lambda}^I} \nc{\klIt}{k_{\lambda}^{\tilde{I}/I}}

\nc{\Rg}{\operatorname{Rg}} \nc{\essinf}{\operatorname{essinf}}
\rnc{\Re}{\operatorname{Re}} 
\nc{\Lin}{\operatorname{Lin}} \nc{\const}{\operatorname{const}}
\nc{\dist}{\operatorname{dist}} \nc{\dpst}{\displaystyle}

\nc{\cl}{c_{\lambda}}
\nc{\dl}{d_{\lambda}}
\nc{\el}{e_{\lambda}}

\nc{\NM}[1]{\|#1\|_{\mathcal{M}(A)}}
\nc{\NH}[1]{\|#1\|_{\mathcal{H}(A)}}
\nc{\NMp}[2]{\|#1\|_{\mathcal{M}(#2)}}
\nc{\NHp}[2]{\|#1\|_{\mathcal{H}(#2)}}

\nc{\Mph}{{\mathcal{M}(\vp)}}
\nc{\Hph}{{\mathcal{H}(\vp)}}

\nc{\Mp}{\mathcal{M}}
\nc{\Hp}{\mathcal{H}}

\title[An extremal problem related to negative refraction]{An extremal problem related to negative refraction}

\author{Kristian Seip}
\address{Department of Mathematical Sciences, Norwegian University of
Science and Technology (NTNU), NO-7491 Trondheim, Norway}
\email{seip@math.ntnu.no}

\author{Johannes Skaar}
\address{Department of Electronics and Telecommunications, Norwegian University of
Science and Technology (NTNU), NO-7491 Trondheim, Norway}
\email{johannes.skaar@iet.ntnu.no}

\thanks{The first author is
supported by the Research Council of Norway grant 160192/V30.}


\begin{document}

\begin{abstract}
We solve an extremal problem that arises in the study of the
refractive indices of passive metamaterials. The problem concerns
Hermitian functions in $H^2$ of the upper half-plane, i.e., $H^2$
functions satisfying $f(-x)=\overline{f(x)}$. An additional
requirement is that the imaginary part of $f$ be nonnegative for
nonnegative arguments. We parameterize the class of such functions
whose real part is constant on an interval, and solve the problem
of minimizing the imaginary part on the interval on which the
function's real part takes a given constant value.
\end{abstract}

\maketitle

\section{Introduction}

We consider in this note an extremal problem that arose in
investigations of certain electromagnetic parameters of artificial
materials (metamaterials). The physical interpretation of our
solution in terms of bounds for refractive indices and theoretical
limitations for the design of metamaterials is described elsewhere
\cite{SS}; the purpose of the present work is to give an account
of the underlying mathematical problem, which seems to be of some
independent interest.

We will be dealing with Hardy spaces $H^p$ of the upper half-plane
$\{z=x+iy:\ y>0\}$. We are primarily interested in $H^2$, but it
is convenient to have at our disposal the whole range of spaces
corresponding to $0<p\le \infty$. For $0<p<\infty$,  $H^p$
consists of those analytic functions $f$ in the upper half-plane
for which
\[ \|f\|_p^p=\sup_{y>0} \int_{-\infty}^\infty |f(x+iy)|^p dx <\infty;\]
$H^\infty$ is the space of bounded analytic functions. A function
$f$ in $H^p$ has a nontangential boundary limit at almost every
point of the real axis, and the corresponding limit function, also
denoted $f$, is in $L^p=L^p(\R)$. Indeed, the $L^p$ norm of the
boundary limit function coincides with the $H^p$ norm introduced
above. Thus we may view $H^p$ as a subspace of $L^p$. We refer to
\cite{Ga} for these and other basic facts about $H^p$, as well as
the twin theory of $H^p$ of the unit disk. (We will make a
reference to the disk setting at one point.)

The Hilbert space $H^2$ is the image of $L^2(\R^+)$ under the
Fourier transform. In practice, it is quite common that one
considers functions in $H^2$ that are Fourier transforms of
real-valued functions in $L^2$. This leads to the following
symmetry condition: $f(-x)=\overline{f(x)}$. Functions $f$
satisfying this condition will be referred to as Hermitian
functions. Thus Hermitian functions have even real parts and odd
imaginary parts.

The Hilbert transform of a function $u$ in $L^p$ ($1\le p<\infty$)
is defined as
\[ \tilde{u}(x)=\text{p.v.}\ \frac{1}{\pi} \int_{-\infty}^\infty
\frac{u(t)}{x-t} dt. \] It acts boundedly on $L^p$ for
$1<p<\infty$ and isometrically on $L^2$. If $u$ is a real-valued
function in $L^p$ for $1<p<\infty$, then $u+i\tilde{u}$ is in
$H^p$, and so the role of the Hilbert transform is to link the
real and imaginary parts of functions in $H^p$. We will only work
with Hermitian functions, and we will be interested in computing
real parts from imaginary parts. For this reason, it will be
convenient for us to consider the following Hilbert operator:
\[\mathcal{H}v(x)=\text{p.v.}\ \frac{1}{\pi}\int_{0}^\infty
v(t)\left(\frac{1}{t-x}+\frac{1}{t+x}\right) dt,\] acting on
functions in $L^p(\R^+)$. Provided $1<p<\infty$, the function
$\mathcal Hv+iv$ will then be in $H^p$, with the presumption that
$v$ is an odd function.

A natural type of problem in $H^2$ is that of approximating a
given Hermitian function supported on two symmetric intervals.
This means that one specifies a desired behavior in a certain
frequency band and attempts to find, in some appropriate sense, an
optimal approximation in $H^2$. Without further restrictions, such
a problem makes little sense, because it is easy to see that
approximations can be made with arbitrary precision in $L^2$ norm.
It may be reasonable to prescribe bounds for the norm of the
approximating function; see for instance the work of M. G.
Kre\u{\i}n and P. Ya. Nudel'man \cite{KrNu1}, \cite{KrNu2} for
interesting results along such lines. In the present note, we take
a different route. We shall require the imaginary part of the
function to be nonnegative for nonnegative
arguments.\footnote{This reflects the passivity condition for our
electromagnetic medium.} We are interested in a specific problem
of this general kind; it turns out to have an explicit and
remarkably simple solution.

\section{Results}

We turn to the statement of the problem. For a finite interval
$I=[a,b]$ ($0<a<b$) and every real number $\alpha$ we define the
family of functions
\[ K_\alpha(I)=\{v\in L^2(\R^+):\ v(t)\ge 0 \
\text{for}  \  t>0, \mathcal{H} v(t)=\alpha \ \text{for} \ t\in
I\}.
\]
(Here and elsewhere we suppress the obvious ``almost everywhere''
provisions needed when considering pointwise restrictions.) We
think of functions in $K_\alpha(I)$, or more generally functions
in $L^2(\R^+)$, as the imaginary parts of Hermitian functions, and
we view them therefore as odd functions on $\R$.

Our purpose is to give a parametrization of $K_{\alpha}(I)$ and to
solve the extremal problem
\[ \lambda=\inf_{v\in K_{-1}(I)}\| \chi_I v\|_\infty, \]
where $\chi_I$ denotes the characteristic function of $I$. We will
show that the extremal problem has the following explicit
solution:
\[ \lambda=\frac{b^2-a^2}{2ab}. \]
We note that the quantity on the right is invariant under
dilations $sI=[s a, s b]$. This is as it should be since $v(t)$ is
in $K_{-1}(I)$ if and only if $v(t/s)$ is in $K_{-1}(sI)$ for
$s>0$. It will become clear that the extremal value $\lambda$ is
not attained by any function in $K_{-1}(I)$. We will also see that
the problem is insensitive to which $L^p$ norm we choose to
minimize.

Clearly, the corresponding extremal problem for $K_\alpha(I)$ has
solution $|\alpha| \lambda$ when $\alpha < 0$. However, if $\alpha
\ge 0$, the extremal problem is uninteresting and has solution
$0$. Thus the sign in the relation $\mathcal{H} v(t)=\alpha$ for
$t\in I$ matters in a decisive way.\footnote{The case $\alpha<-1$
corresponds to the interesting physical phenomenon of negative
refraction, which has received considerable attention in recent
years. Artificial, negatively refracting materials, called
metamaterials, have been realized in the microwave range
\cite{Sm}, building on previous theoretical ideas
\cite{Ve,Pe1,Pe2}. As explained in \cite{SS}, the solution to our
extremal problem provides a bound for the loss of negatively
refracting materials when the real part of the refractive index is
constant in a finite bandwidth.}

The following lemma is basic for our parametrization of
$K_\alpha(I)$.

\begin{lemma}
A real-valued function in $L^2(\R^+\setminus I)$ is the
restriction to $\R^+\setminus I$ of at most one real-valued
function $v$ in $L^2(\R^+)$ such that $\mathcal{H}v(t)$ is
constant on $I$.
\end{lemma}

\begin{proof}
We assume two real-valued functions $v_1$ and $v_2$ in $L^2(\R^+)$
coincide off $I$ and are such that both $u_1(t)=\mathcal{H}v_1(t)$
and $u_2(t)=\mathcal{H}v_2(t)$ are constant on $I$. If we set
$c=u_2(t)-u_1(t)$ for $t$ in $I$, then the function
$h=[(u_1-u_2+i(v_1-v_2)+c]^2$ will be real for real arguments. A
change of variables argument shows that then $h(i(1+z)/(1-z))$
belongs to $H^1$ of the unit disk. But a function in $H^1$ can be
real only if it is a constant. Clearly, $h$ can be a constant only
if $u_1-u_2+i(v_1-v_2)=0$.
\end{proof}

We note that the assumption that $v$ is in $L^2$ is essential for
this lemma; the proof would break down if we assumed, say, that
$v$ belonged to some $L^p$ for $p<2$.

The following function will play an essential role in what
follows:
\[ \sigma(z)=\frac{1}{\sqrt{z^2-b^2}\sqrt{z^2-a^2}}. \]
This function, which is taken to be positive for real arguments
$x>b$, is analytic in the slit plane
$\C\setminus([-b,-a]\cup[a,b])$. For real arguments $a<|x|<b$ we
define $\sigma(x)$ by extending it continuously from the upper
half-plane. Thus $\sigma(x)$ takes values on the negative
imaginary half-axis when $x$ is in $(a,b)$ and on the positive
imaginary half-axis when $-x$ is in $(a,b)$, and otherwise it is
real for real arguments. The key point, besides the symmetry
$\sigma(-x)=\overline{\sigma(x)}$, is that $\sigma$ provides a
means for switching between real and imaginary when switching off
and on $I$.

The following is our main result. \pagebreak

\begin{theorem*}
A nonnegative function $v$ in $L^2(\R^+)$ is in $K_\alpha(I)$ if
and only if the following three conditions hold:

\begin{equation} \label{one} \int_{\R^+\setminus I} v(t) |\sigma(t)| dt < \infty
\end{equation}
\begin{equation} \label{two} \frac{2}{\pi}\int_{\R^+\setminus I} tv(t) \sigma(t) dt =
\alpha \end{equation} \begin{equation} \label{three}
v(x)=\mathcal{H}\left((1-\chi_I)\sigma v\right)(x)/|\sigma(x)|, \
\ \ x\in I.
\end{equation}

\end{theorem*}

Some remarks are in order before we give the proof of the theorem.

The integrability condition \eqref{one} is merely a slight growth
condition at the endpoints of $I$; we may write it more succinctly
as
\[
\int_0^a\left[v(a-t)+v(b+t)\right]\frac{dt}{\sqrt{t}} <\infty.
\]
This condition ensures that the integral in \eqref{two} and the
Hilbert transform appearing in \eqref{three} are both
well-defined.

At first sight, the theorem may not seem to give an explicit
parametrization of $K_\alpha(I)$. However, the Hilbert transform
appearing in \eqref{three} is given by
\[ \mathcal{H}\left((1-\chi_I)\sigma
v\right)(x)=\frac{1}{\pi}\int_{\R^+\setminus I} v(t)
\sigma(t)\frac{2t}{t^2-x^2}  dt, \] and we observe that the
integrand on the right is nonnegative whenever $v(t)$ is
nonnegative. Hence $v(x)\ge 0$ for $x$ off $I$ implies $v(x)\ge 0$
for $x$ in $I$. This small miracle implies that $K_\alpha(I)$ is
parameterized by those nonnegative functions $v$ in
$L^2(\R^+\setminus I)$ for which \eqref{one} and \eqref{two} hold
and such that
\[ \int_I |\mathcal{H}\left((1-\chi_I)\sigma v\right)(x)|^2 |\sigma(x)|^{-2}
dx <\infty. \] By rephrasing this condition in more explicit terms
(see Lemma~3 below), we arrive at the following corollary.

\begin{corollary}
A nonnegative function $\nu$ in $L^2(\R^+\setminus I)$ has an
extension to a function in some class $K_\alpha(I)$ if and only if
the following condition holds:
\[ \int_0^a
\int_0^a\left[\nu(a-t)\nu(a-\tau)+\nu(b+t)\nu(b+\tau)\right]\frac{|\log(t+\tau)|}{\sqrt{t\tau}}
\ dtd\tau<\infty. \]
\end{corollary}
The difference between \eqref{one} and the condition above is the
logarithmic factor, which means that the condition of the
corollary is only a very slight strengthening of \eqref{one}. It
is clear that for instance boundedness of $v$ near the endpoints
of $I$ is more than enough.

We note that the integrand in \eqref{two} is negative to the left
of $I$ and positive to the right of $I$. This means that if
$\alpha$ is negative, then
\[ |\alpha|\le \frac{2}{\pi}\int_{0}^a tv(t) |\sigma(t)| dt, \]
with equality holding if $v$ vanishes to the right of $I$. It
follows that \[\mathcal{H}\left((1-\chi_I)\sigma v\right)(x)\ge
\frac{1}{\pi}\int_0^a v(t) \sigma(t)\frac{2t}{t^2-x^2} dt \ge
\frac{|\alpha|}{x^2};
\]
we may come as close as we wish to this lower bound by choosing
any suitable $v$ supported on a small set sufficiently close to
$0$. Hence our extremal problem (corresponding to $\alpha=-1$) has
solution
\[ \lambda=\max_{x\in I}\frac{\sqrt{(b^2-x^2)(x^2-a^2)}}{x^2}=
\frac{b^2-a^2}{2ab}, \] as proclaimed above. We also observe that
the same function $1/(x^2|\sigma(x)|)$ would give the infimum for
the $L^p$ norm over $I$ for any other value of $p>0$.

If, on the other hand, $\alpha$ is positive, we have instead
\[ \alpha\le \frac{2}{\pi}\int_{b}^\infty tv(t) |\sigma(t)| dt, \]
with equality holding if $v$ vanishes to the left of $I$. In this
case, arguing in the same fashion as above, we find that we can
get $\mathcal{H}\left((1-\chi_I)\sigma v\right)(x)$ as small as we
please by letting $v$ be supported on a set sufficiently far to
the right of $I$.

\section{Proofs}

We now turn to the proof of the theorem and its corollary. We will
rely on Lemma~1 and two additional lemmas.

\begin{lemma}

For every $t$ in $(0,a)\cup (b,\infty)$ the function \[f_t(x)=
\frac{2t}{t^2-x^2}\left(1-\frac{\sigma(t)}{\sigma(x)}\right)-2t\sigma(t)
\] is in $H^p$ for $p>1/2$, and the following estimates hold
\[ \|f_t\|_1 \le \frac{C_1}{\sqrt{|(a-t)(b-t)|}}, \ \ t<2b,  \]
\[ \int_{2b}^\infty |f_t(x)|^2 dt \le \frac{C_2}{|x|+1}, \]
where the constants $C_1$ and $C_2$ only depend on $a$ and $b$.
\end{lemma}

\begin{proof} It is immediate that $f_t$ belongs to $H^p$ for $p>1/2$ because the
isolated singularities $\pm t$ are removable and
$f_t(z)=O(z^{-2})$ when $z\to\infty$. The norm estimates follow
from elementary calculations. \end{proof}

\begin{lemma}
A nonnegative function $\nu$ in $L^2(\R^+\setminus I)$ satisfies
\[ \int_I |\mathcal{H}((1-\chi_I)\sigma \nu)(x)|^2|\sigma(x)|^{-2}
dx <\infty \] if and only if the following condition holds:
\[ \int_0^a
\int_0^a\left[\nu(a-t)\nu(a-\tau)+\nu(b+t)\nu(b+\tau)\right]\frac{|\log(t+\tau)|}{\sqrt{t\tau}}
\ dtd\tau<\infty. \]
\end{lemma}
\begin{proof}
The necessary and sufficient condition for square-integrability at
the left end-point of $I$ is that
\[ \int_{0}^a\left(\int_0^a \frac{\nu(a-t)}{\sqrt{t}(t+x)} dt\right)^2
x\ dx <\infty. \] By Fubini's theorem, we may interchange the
order of integration so that this condition becomes
\[ \int_0^a
\int_0^a \nu(a-t)\nu(a-\tau)\frac{|\log(t+\tau)|}{\sqrt{t\tau}} \
dtd\tau<\infty. \] Combining this with the corresponding condition
at the right end-point of $I$, we arrive at the condition of the
lemma.
\end{proof}

The theorem is now proved in the following way. We assume first
that we are given a nonnegative function $v$ in $L^2(\R^+)$
satisfying \eqref{one}, \eqref{two}, and \eqref{three}. We claim
that the function
\[f_1(x)=\frac{1}{\pi}\int_{(0,a)\cup(b,2b)} v(t)f_t(x)dt\]
is in $H^1$. Indeed, by Lemma~2 and Fubini's theorem,
\[ \| f_1 \|_1\le \frac{C_1}{\pi} \int_{(0,a)\cup(b,2b)}
\frac{v(t)}{\sqrt{|(a-t)(b-t)|}} dt, \] and the integral on the
right is bounded thanks to \eqref{one}. On the other hand,
\[f_2(x)=\frac{1}{\pi}\int_{2b}^\infty v(t)f_t(x) dt\]
is in $H^p$ for $p>2$, because by the Cauchy--Schwarz inequality
and Lemma~2 we have
\[ \int_{-\infty}^\infty |f_2(x)|^p dx \le \frac{C_2^{p/2}}{\pi^p} \|v\|_2^p
\int_{-\infty}^\infty \frac{1}{(|x|+1)^{p/2}}dx. \] The imaginary
part of $f_1+f_2$ is supported by $I$ and equals $-v$ there, in
view of \eqref{three}. Since $v$ is assumed to be in $L^2$, it
follows that $f_1+f_2$ is in fact in $H^2$.

We set $f=\mathcal{H}((1-\chi_I)v)+i(1-\chi_I)v$, which is a
function in $H^2$. We observe that the imaginary part of
$f-(f_1+f_2)$ equals $v$ and that its real part equals $\alpha$ on
$I$, when taking into account \eqref{two}. So we have proved that
the given $v$ is indeed in $K_\alpha(I)$.

We now prove the necessity of the three conditions of the theorem.
So assume we are given a nonnegative function $v$ in $L^2(\R^+)$
belonging to some class $K_\alpha(I)$. Setting
$I_\eps=[a+\eps,b-\eps]$, we see that $v$ also belongs to
$K_\alpha(I_\eps)$ whenever $0<\eps<(b-a)/2$. But since the
Hilbert transform of $v$ is constant near the endpoints of
$I_\eps$, it follows that
\[
\int_{\R^+\setminus I_\eps} v(t) |\sigma_\eps(t)| dt < \infty,
\]
where now
\[
\sigma_\eps(t)=\frac{1}{\sqrt{t^2-(b-\eps)^2}\sqrt{t^2-(a+\eps)^2}}.
\]
We claim that this means that \begin{equation} \label{infty} v(x)=
\mathcal{H}\left((1-\chi_{I_\eps})\sigma_\eps
v\right)(x)/|\sigma_\eps(x)|,\end{equation} provided $x$ is in
$(a+\eps,b-\eps)$. Indeed, by Lemma~1, it is enough to verify that
\[ (1-\chi_{I_\eps}(x))v(x)+\chi_{I_\eps}(x)
\mathcal{H}\left((1-\chi_{I_\eps})\sigma_\eps
v\right)(x)/|\sigma_\eps(x)|\] is in $K_\alpha(I_\eps)$ for some
$\alpha$. Since, in view of Lemma~3, the function on the
right-hand side of \eqref{infty} is square-integrable on $I_\eps$,
the claim follows by repeating the argument in the first part of
the proof.

We may view $\mathcal{H}\left((1-\chi_{I_\eps})\sigma_\eps
v\right)(x)$ as the $L^1$ norm of the function
\[
h_{x,\eps}(t)=\frac{1}{\pi}\left((1-\chi_{I_\eps}(t))\sigma_\eps(t)
v(t)\right)\frac{2t}{t^2-x^2}. \] Then \eqref{infty} says that
$\|h_{x,\eps}\|_1 \to v(x)|\sigma(x)|$ when $\eps\to 0$. Since we
also have that
\[ h_{x,\eps}(t)\to \frac{1}{\pi}\left((1-\chi_{I}(t))\sigma(t)
v(t)\right)\frac{2t}{t^2-x^2} \] for every $t$, we obtain
\[
 v(x)=\mathcal{H}\left((1-\chi_{I})\sigma v\right)(x)/|\sigma(x)|\]for
every $x$ in $(a,b)$. By a similar argument, we find that
\[ \alpha=\lim_{\eps\to 0} \frac{2}{\pi}\int_{\R^+\setminus I_\eps}
tv(t) \sigma_\eps(t) dt=\frac{2}{\pi}\int_{\R^+\setminus I} tv(t)
\sigma(t) dt. \] The necessity of \eqref{one} has already been
observed; without it we would reach the contradictory conclusion
that $v(x)=\infty$ for almost every $x\in (a,b)$.

We finally note that the corollary is an immediate consequence of
the theorem and lemmas 1 and 3.

\end{document}